\documentclass[12pt]{article}

\usepackage[a4paper,margin=1in]{geometry}
\usepackage{amsmath,amssymb,amsthm,mathtools}
\usepackage{mathrsfs}
\usepackage{enumitem}
\usepackage{microtype}
\usepackage{hyperref}
\usepackage[nameinlink,capitalise]{cleveref}

\hypersetup{
  colorlinks=true,
  linkcolor=blue,
  citecolor=blue,
  urlcolor=blue
}

\newtheorem{definition}{Definition}[section]
\newtheorem{lemma}[definition]{Lemma}
\newtheorem{proposition}[definition]{Proposition}

\newtheorem{corollary}[definition]{Corollary}
\newtheorem{remark}[definition]{Remark}

\newcommand{\ord}{\operatorname{ord}}
\newcommand{\Hom}{\operatorname{Hom}}

\newcommand{\id}{\operatorname{id}}

\title{\Large 
An Exposition of Gerstenhaber's Proposition 1: Elimination of Coboundaries in Formal Deformations}
 \usepackage{authblk}

\author{Rinkila Bhutia}
\author{Zuiner Jamir}

\affil{Department of Mathematics, Sikkim University, Gangtok, 737102, Sikkim, India}
\date{}

\begin{document}
\maketitle

\begin{abstract}
Gerstenhaber's foundational work on formal deformations of associative algebras
relates the first non-trivial term of a deformation to the second Hochschild
cohomology group.  In the proof of Proposition~1 of his 1964 paper, a leading
Hochschild $2$-cocycle which is a coboundary is successively removed by formal
changes of variables.  If this elimination terminates after finitely many
steps, no limiting argument is required, but the possibility that coboundaries continue to be eliminated indefinitely is left largely implicit in the original treatment. This expository note makes that limiting argument fully explicit for deformations over $K[[t]]$ of a fixed associative $K$-algebra $A$. We use the $t$-adic topology on the relevant spaces of formal
power series and formal endomorphisms, establish completeness
and the continuity of composition, inversion and the action
on formal multiplications, and show that an infinite sequence
of successive changes of variables assembles into a single,
well-defined formal equivalence. This yields a complete, self-contained proof of the classical result and its rigidity consequence. 
\end{abstract}


\tableofcontents

\section{Introduction}

The formal deformation theory of associative algebras developed by
Gerstenhaber provides a systematic way of studying how the multiplication
of an algebra can be varied while preserving associativity. One of the
central features of the theory is the appearance of Hochschild cohomology.
For an associative algebra $A$, the first-order term of a formal deformation
is a Hochschild $2$-cocycle, and changing the deformation by a formal
change of variables changes this cocycle by a Hochschild coboundary.
Thus its cohomology class belongs naturally to $H^2(A,A)$
\cite{Hochschild1945,Gerstenhaber1964}.

Gerstenhaber's 1963 paper on the cohomology structure of associative
rings established important cohomological structures on Hochschild
cohomology. He subsequently developed the deformation theory of rings
and algebras in his 1964 paper and continued this work in three further
papers \cite{Gerstenhaber1963,Gerstenhaber1964,Gerstenhaber1966,
Gerstenhaber1968,Gerstenhaber1974}.

The purpose of this paper is to give a detailed exposition of the argument
behind Proposition~1 of Gerstenhaber's paper \cite{Gerstenhaber1964}.  In the
notation used below, a formal one-parameter deformation has the form
\[
  \mu_t=\mu+tF_1+t^2F_2+\cdots,
\]
where $\mu$ is the original multiplication.  Gerstenhaber's proposition says
that,  replacing the deformation by an equivalent one, it can be written in the form
%(not sure if this is the right terminology), it can be written: uses equivalent in pg 65 and also in proposition 1
\[
  \mu_t=\mu+t^nF_n+t^{n+1}F_{n+1}+\cdots,
\]
where
\[
  F_n\in Z^2(A,A),
  \qquad [F_n]\neq 0\in H^2(A,A).
\]
In particular, the first nonvanishing coefficient (\emph{infinitesimal}) can be chosen not to be a
coboundary.

The following is the idea behind this assertion.  Suppose the infinitesimal is $F_n$ and that
\[
  F_n=\delta q
\]
is a coboundary.  A formal change of variables of the form
\[
  \phi_t=\id_A+t^nq
\]
%(with the sign chosen according to the convention for the Hochschild differential)
removes the coefficient of $t^n$.  The first nonzero term is
therefore pushed to a strictly higher power of $t$.

If this procedure stops after finitely many steps, the argument is immediate:
one eventually reaches a leading cocycle which is not a coboundary.  The
subtle point is the possibility that the procedure continues indefinitely:
\[
  \mu^{(0)}\sim\mu^{(1)}\sim\mu^{(2)}\sim\cdots,
\]
with a coboundary eliminated at every stage.  In that situation there is no
finite stage at which one can simply declare the process finished.  One must
show that the infinitely many elementary equivalences combine to a single
well-defined formal equivalence.

This is the point on which the present exposition concentrates.  We show that
the successive changes of variables begin at strictly increasing powers of
$t$.  Consequently their partial compositions form a Cauchy sequence in the
$t$-adic topology.  Completeness then gives a limiting formal automorphism %(i guess we can change this term also),
and a parallel argument gives its inverse.  Continuity of composition and of
the action on the multiplication allows us to pass to the limit.

The organization of the paper reflects this purpose. Section 2 recalls the formal deformation and equivalence framework and proves the elementary elimination lemma for a coboundary. Section 3 develops the \(t\)-adic framework and establishes completeness. Section 4 studies formal endomorphisms, continuity of composition and inversion, the action on formal multiplications, and infinite composition of elementary transformations. Section 5 applies these results to the infinite elimination process and gives the resulting proof of Gerstenhaber’s Proposition 1, followed by the rigidity consequence.
%Section~9 comments explicitly on the relationship between the concise original proof and the expanded argument given here.

\section{Formal deformations and equivalence}\label{formal deformation}

Let $A$ be an associative algebra over a field $K$, with multiplication
\[
  \mu\colon A\otimes A\longrightarrow A,
  \qquad \mu(a,b)=ab.
\]
We write $C^n(A,A)=\Hom_K(A^{\otimes n},A)$ for the Hochschild $n$-cochains.

\subsection{Formal one-parameter deformations}

\begin{definition}
A one-parameter formal deformation of the multiplication $\mu$ is a formal
power series
\[
  \mu_t=\mu+tF_1+t^2F_2+\cdots,
  \qquad F_i\in C^2(A,A),
\]
which is associative over $A[[t]]$.
\end{definition}

Associativity means
\[
  \mu_t(\mu_t(a,b),c)=\mu_t(a,\mu_t(b,c))
\]
for all $a,b,c\in A$, as an identity of formal power series in $t$.

We use the Hochschild differential in the convention
\[
\begin{aligned}
(\delta q)(a,b)
 &=aq(b)-q(ab)+q(a)b,
 && q\in C^1(A,A),\\
(\delta F)(a,b,c)
 &=aF(b,c)-F(ab,c)+F(a,bc)-F(a,b)c,
 && F\in C^2(A,A).
\end{aligned}
\]

Comparing the coefficient of $t$ in the associativity equation gives
\[
  \delta F_1=0.
\]
Thus
\[
  F_1\in Z^2(A,A).
\]
More generally, suppose
\[
  F_1=\cdots=F_{n-1}=0
\]
and $F_n$ is the first nonzero coefficient, we call $F_n$ the \emph{infinitesimal}.  Comparing the coefficient of
$t^n$ in the associativity equation again gives
\[
  \delta F_n=0.
\]
Hence the infinitesimal is always a Hochschild $2$-cocycle.

\subsection{Equivalence of formal deformations}

\begin{definition}
Two formal deformations $\mu_t$ and $\nu_t$ are called equivalent if there is a
formal $K[[t]]$-linear automorphism
\[
  \phi_t=\id_A+t\phi_1+t^2\phi_2+\cdots
\]
of $A[[t]]$ such that
\[
  \nu_t
  =\phi_t\,\mu_t\,(\phi_t^{-1}\otimes\phi_t^{-1}).
\]
\end{definition}

Since the constant term of $\phi_t$ is $\id_A$, the formal inverse
\[
\phi_t^{-1}=\id_A+t\psi_1+t^2\psi_2+\cdots
\]
exists uniquely.\\
%(will have to say what $v_t$ is)
Suppose
\[
\mu_t=\mu+tF_1+t^2F_2+\cdots
\]
and
\[
\nu_t=\mu+tG_1+t^2G_2+\cdots.
\]
Expanding the equivalence relation to first order gives
\[
G_1=F_1-\delta\phi_1.
\]
Consequently,
\[
[G_1]=[F_1]\in H^2(A,A).
\]
Thus the Hochschild cohomology class of the first order term is an invariant of
the equivalence class of the deformation.

If $F_1=0$, the deformation may nevertheless be nontrivial, with its first
nonzero term occurring at some higher order. The following lemma shows that a
coboundary occurring as the leading term can be eliminated by a formal change
of variables.

\begin{lemma}[Elimination lemma]\label{elimination lemma}
Suppose that
\[
\mu_t=\mu+t^nF_n+t^{n+1}F_{n+1}+\cdots
\]
and that
\[
F_n=\delta q
\]
for some $q\in C^1(A,A)$. Define
\[
\phi_t=\id_A+t^nq.
\]
Then the equivalent deformation
\[
\mu'_t
=\phi_t\mu_t(\phi_t^{-1}\otimes\phi_t^{-1})
\]
has no term of degree $n$; that is,
\[
\mu'_t=\mu+O(t^{n+1}).
\]
\end{lemma}

\begin{proof}
 Write the formal inverse of $\phi_t$ in the form
\[
\phi_t^{-1}
=\id_A+t^1\psi_1+t^2\psi_2+\cdots
\]
Using
\[
\phi_t\phi_t^{-1}=\id_A
\]
and comparing the coefficient of $t^i$, we obtain
\[\psi_i=0~for~i<n~and~
q+\psi_n=0.
\]
Hence
\[
\phi_t^{-1}=\id_A-t^nq+O(t^{n+1}).
\]
Substituting into the transformed multiplication gives
\[
\mu'_t
=(\id_A+t^nq)
\bigl(\mu+t^nF_n+O(t^{n+1})\bigr)
\circ
\bigl((\id_A-t^nq+O(t^{n+1}))
\otimes
(\id_A-t^nq+O(t^{n+1}))\bigr).
\]
Only the terms of degree $n$ are relevant. For $a,b\in A$, the coefficient of
$t^n$ is
\[
F_n(a,b)+q(ab)-q(a)b-aq(b).
\]
By the definition of the Hochschild differential,
\[
(\delta q)(a,b)
=aq(b)-q(ab)+q(a)b,
\]
and hence
\[
q(ab)-q(a)b-aq(b)=-\delta q(a,b).
\]
Therefore the coefficient of $t^n$ in $\mu'_t$ is
\[
F_n-\delta q=0.
\]
Thus
\[
\mu'_t=\mu+O(t^{n+1}),
\]
as required.
\end{proof}
The elimination lemma shows that a suitable formal change of variables can eliminate a coboundary occurring at a prescribed order of a deformation. To apply this idea successively, we need to keep track of the order of formal
power series and of the effect of successive changes of variables on their higher-order terms. This leads naturally to the $t$-adic viewpoint. We therefore introduce the necessary
$t$-adic framework before returning to Gerstenhaber's Proposition~1.

\section{The $t$-adic framework}

Let $V$ be a $K$-vector space.  We use the formal power series space
\[
  V[[t]]
  =\left\{\sum_{k\geq0}v_kt^k:v_k\in V\right\}.
\]

% We equip $V[[t]]$ with the $t$-adic topology, that is, the topology
% defined by the descending sequence of subspaces
% \[
%   V[[t]]\supset tV[[t]]\supset t^2V[[t]]
%   \supset\cdots.
% \]
% This is the special case of the general $I$-adic topology corresponding
% to the ideal $(t)$; see, for example, \cite{Matsumura1989}.
% Equivalently, the $t$-adic topology may be described by the metric
% \[
%   d(x,y)=2^{{ord}(x-y)},
% \]
% with the convention that $d(f,f)=0$.
% %\subsection{Order and the $t$-adic metric}

For a nonzero series
\[
  x=\sum_{k\geq0}v_kt^k,
\]
define
\[
  \ord(x)=\min\{k\geq0:v_k\neq0\},
\]
and set $\ord(0)=+\infty$.  Define
\[
  d(x,y)=2^{-\ord(x-y)}.
\]

\begin{lemma}
The function $d$ is an ultrametric on $V[[t]]$.
\end{lemma}

\begin{proof}
The only nontrivial point is the strong triangle inequality. 
For $x,y,z\in V[[t]]$, we have

$$
x-z=(x-y)+(y-z).
$$

Hence, using

\begin{equation}
    \ord(u+v)\geq \min\{\ord(u),\ord(v)\}, \label{one}
\end{equation}

we obtain

$$
\ord(x-z)
\geq
\min\{\ord(x-y),\ord(y-z)\}.
$$

Now, if $f$ is a decreasing function, then

$$
f(\min\{a,b\})=\max\{f(a),f(b)\}.
$$

Taking $f(r)=2^{-r}$, which is decreasing, we therefore have

$$
\begin{aligned}
d(x,z)
&=2^{-\ord(x-z)}\\
&\leq
2^{-\min\{\ord(x-y),\ord(y-z)\}}\\
&=
\max\left\{
2^{-\ord(x-y)},2^{-\ord(y-z)}
\right\}\\
&=
\max\{d(x,y),d(y,z)\}.
\end{aligned}
$$
The remaining metric properties are immediate.
Thus $d$ satisfies the strong triangle inequality and is therefore an
ultrametric.

\end{proof}

\subsection{Completeness}

\begin{proposition}
The metric space $(V[[t]],d)$ is complete.
\end{proposition}

\begin{proof}
Let $(x_i)$ be a Cauchy sequence in $V[[t]]$, where

$$
x_i=\sum_{k\geq0}x_{i,k}t^k.
$$

For every $N\geq0$, since $(x_i)$ is Cauchy, there exists $I_N$ such that
for all $i,j\geq I_N$,

$$
d(x_i,x_j)<2^{-N}.
$$

By the definition of $d$, this implies

$$
2^{-\ord(x_i-x_j)}<2^{-N},
$$

and hence

$$
\ord(x_i-x_j)>N.
$$

Since the order is integer-valued, we have

$$
\ord(x_i-x_j)\geq N+1.
$$

Therefore,

$$
x_i\equiv x_j\pmod{t^{N+1}}.
$$

In particular, the coefficients of $x_i$ and $x_j$ agree in all degrees
$0,1,\ldots,N$. Hence, for every fixed $k$, the coefficient $x_{i,k}$ is
eventually constant. Let its eventual value be $v_k$, and define

$$
x=\sum_{k\geq0}v_kt^k\in V[[t]].
$$

For any $N\geq0$, the coefficients of $x_i$ and $x$ agree in degrees
$0,1,\ldots,N$ once $i$ is sufficiently large. Thus

$$
\ord(x_i-x)\geq N+1,
$$

and consequently

$$
d(x_i,x)
\leq 2^{-(N+1)}.
$$

Since $N$ is arbitrary, it follows that

$$
x_i\longrightarrow x.
$$

Therefore $(V[[t]],d)$ is complete.
\end{proof}

\section{Formal endomorphisms and infinite composition}

The preceding completeness result will now be applied to formal power
series whose coefficients are linear endomorphisms.

Put
\[
  \mathcal E=\Hom_K(A,A)[[t]].
\]

Thus an element of $\mathcal E$ is a formal series
$
\Phi=\sum_{r\geq0}\Phi_r t^r,~
\Phi_r\in\operatorname{Hom}_K(A,A)$.
%Equivalently, such a series may be regarded as a formal $K[[t]]$-linear endomorphism of $A[[t]]$.

If $
  p=\sum_{r\geq0}p_rt^r,
  \text{ and }
  q=\sum_{s\geq0}q_st^s,
$
then define composition on $\mathcal E$ as
\[
  p\circ q
  =\sum_{k\geq0}
    \left(\sum_{r+s=k}p_r\circ q_s\right)t^k.
\]

The metric introduced in the previous section gives a natural topology on
$\mathcal E$.

We first record the basic estimate governing composition.

\begin{lemma}
For $p,q\in\mathcal E$,
\[
  \ord(p\circ q)
  \geq \ord(p)+\ord(q).
\]
\end{lemma}

\begin{proof}
If $k<\ord(p)+\ord(q)$ and $r+s=k$, then either
$r<\ord(p)$ or $s<\ord(q)$.  Thus either $p_r=0$ or $q_s=0$.  Hence every
coefficient of degree less than $\ord(p)+\ord(q)$ vanishes.
\end{proof}

We will also repeatedly use the elementary inequality \ref{one}
\[
\ord(f+g)\geq\min\{\ord(f),\ord(g)\}.
\]
In particular, if both $\ord(f_i)$ and $\ord(g_i)$ tend to $+\infty$,
then
\[
\ord(f_i+g_i)\longrightarrow+\infty.
\]
Further by the coefficientwise definition of the tensor product, we have

$$
{ord}(\psi\otimes\psi')
\geq
{ord}(\psi)+{ord}(\psi').
$$

\subsection{Continuity of composition}

\begin{proposition}
If $p_i\to p$ and $q_i\to q$ in $\mathcal E$, then
\[
  p_i\circ q_i\longrightarrow p\circ q.
\]
\end{proposition}

\begin{proof}
We have
\[
p_i\circ q_i-p\circ q
=
(p_i-p)\circ q_i
+
p\circ(q_i-q).
\]

Since
\[
p_i\to p,
\]
we have
\[
\ord(p_i-p)\longrightarrow+\infty.
\]
Also, since $q_i\in\mathcal E=\operatorname{Hom}_K(A,A)[[t]]$, we have
\[
\ord(q_i)\geq0.
\]
Therefore, by the order estimate for composition,
\[
\ord\bigl((p_i-p)\circ q_i\bigr)
\geq
\ord(p_i-p)+\ord(q_i)
\geq
\ord(p_i-p),
\]
and hence
\[
\ord\bigl((p_i-p)\circ q_i\bigr)
\longrightarrow+\infty.
\]

Similarly,
\[
\ord(q_i-q)\longrightarrow+\infty,
\]
and since $\ord(p)\geq0$,
\[
\ord\bigl(p\circ(q_i-q)\bigr)
\geq
\ord(p)+\ord(q_i-q)
\geq
\ord(q_i-q).
\]
Thus
\[
\ord\bigl(p\circ(q_i-q)\bigr)
\longrightarrow+\infty.
\]

Using
\[
\ord(f+g)
\geq
\min\{\ord(f),\ord(g)\},
\]
we obtain
\[
\ord\bigl(p_i\circ q_i-p\circ q\bigr)
\longrightarrow+\infty.
\]
Therefore
\[
p_i\circ q_i\longrightarrow p\circ q.
\]
\end{proof}

% \begin{proposition}
% If $p_i\to p$ and $q_i\to q$ in $\mathcal E$, then
% \[
%   p_i\circ q_i\longrightarrow p\circ q.
% \]
% \end{proposition}

% \begin{proof}
% Write
% \[
%   p_i\circ q_i-p\circ q
%   =(p_i-p)\circ q_i+p\circ(q_i-q).
% \]
% Since $q_i\to q$, the orders $\ord(q_i)$ are bounded below by a fixed integer
% for all sufficiently large $i$.  The order (order estimate?))therefore gives
% \[
%   \ord\bigl((p_i-p)\circ q_i\bigr)\longrightarrow\infty.
% \]
% Similarly,
% \[
%   \ord\bigl(p\circ(q_i-q)\bigr)\longrightarrow\infty.
% \]
% Hence the difference tends to zero.\\(here the fact that ord(f+g)$\geq$ min(ord(f),ord(g)) is being used)
% \end{proof}

\subsection{Continuity of inversion}

\begin{proposition}
Suppose
\[
\Phi_i\longrightarrow\Phi
\]
in $\mathcal E$, and suppose that each $\Phi_i$ has the form
\[
\Phi_i=\id_A+O(t).
\]
Then
\[
\Phi_i^{-1}\longrightarrow\Phi^{-1}.
\]
\end{proposition}

\begin{proof}
Write
\[
\Phi=\id_A+\sum_{m\geq1}\Phi_m t^m
\]
and
\[
\Phi^{-1}
=
\id_A+\sum_{m\geq1}\Psi_m t^m.
\]
Since
\[
\Phi\circ\Phi^{-1}=\id_A,
\]
we have
\[
\left(
\id_A+\sum_{m\geq1}\Phi_m t^m
\right)
\circ
\left(
\id_A+\sum_{m\geq1}\Psi_m t^m
\right)
=
\id_A.
\]

Comparing coefficients of $t^k$, for $k\geq1$, gives
\[
\Psi_k+\sum_{r=1}^{k}\Phi_r\circ\Psi_{k-r}=0.
\]
Hence
\[
\Psi_k
=
-\sum_{r=1}^{k}\Phi_r\circ\Psi_{k-r}.
\]
In particular,
\[
\Psi_1=-\Phi_1,
\]
and, more generally, $\Psi_k$ is determined recursively by
\[
\Phi_1,\Phi_2,\ldots,\Phi_k.
\]

Consequently, if
\[
\Phi=\id_A+\sum_{m\geq1}\Phi_m t^m
\]
and
\[
\widetilde{\Phi}
=
\id_A+\sum_{m\geq1}\widetilde{\Phi}_m t^m
\]
satisfy
\[
\Phi_m=\widetilde{\Phi}_m
\qquad\text{for }1\leq m\leq N,
\]
then their inverses
\[
\Phi^{-1}
=
\id_A+\sum_{m\geq1}\Psi_m t^m
\]
and
\[
\widetilde{\Phi}^{-1}
=
\id_A+\sum_{m\geq1}\widetilde{\Psi}_m t^m
\]
satisfy
\[
\Psi_m=\widetilde{\Psi}_m
\qquad\text{for }1\leq m\leq N.
\]

Indeed, this follows by induction on $m$. For $m=1$,
\[
\Psi_1=-\Phi_1
=
-\widetilde{\Phi}_1
=
\widetilde{\Psi}_1.
\]
Suppose that
\[
\Psi_j=\widetilde{\Psi}_j
\qquad\text{for }1\leq j<m.
\]
Then
\[
\begin{aligned}
\Psi_m
&=
-\sum_{r=1}^{m}\Phi_r\circ\Psi_{m-r}\\
&=
-\sum_{r=1}^{m}\widetilde{\Phi}_r
\circ\widetilde{\Psi}_{m-r}\\
&=
\widetilde{\Psi}_m,
\end{aligned}
\]
because $\Phi_r=\widetilde{\Phi}_r$ for $r\leq m$ and, by the induction
hypothesis, $\Psi_{m-r}=\widetilde{\Psi}_{m-r}$ for $m-r<m$.

Thus, if
\[
\Phi\equiv\widetilde{\Phi}\pmod{t^{N+1}},
\]
then
\[
\Phi^{-1}\equiv\widetilde{\Phi}^{-1}\pmod{t^{N+1}}.
\]

We now apply this to the sequence $(\Phi_i)$. Since
\[
\Phi_i\longrightarrow\Phi,
\]
for every $N\geq0$ there exists $I_N$ such that, for all $i\geq I_N$,
\[
d(\Phi_i,\Phi)<2^{-N}.
\]
Hence
\[
\Phi_i\equiv\Phi\pmod{t^{N+1}}.
\]
By the coefficientwise argument above,
\[
\Phi_i^{-1}
\equiv
\Phi^{-1}
\pmod{t^{N+1}}.
\]
Therefore
\[
\ord(\Phi_i^{-1}-\Phi^{-1})\geq N+1,
\]
and consequently
\[
d(\Phi_i^{-1},\Phi^{-1})
\leq
2^{-(N+1)}.
\]

Since $N$ is arbitrary, it follows that
\[
\Phi_i^{-1}\longrightarrow\Phi^{-1}.
\]
\end{proof}
% \subsection{Formal inverses}

% \begin{lemma}
% If
% \[
%   \Phi=\id_A+R,
%   \qquad \ord(R)\geq1,
% \]
% then $\Phi$ is invertible in $\mathcal E$, with
% \[
%   \Phi^{-1}
%   =\id_A-R+R^{\circ2}-R^{\circ3}+\cdots,
% \]
% where $R^{\circ k}$ denotes the $k$-fold composition of $R$ with itself.
% \end{lemma}

% \begin{proof}
% By the composition-order inequality,
% \[
%   \ord(R^{\circ k})\geq k.
% \](this is because the order is at least 1)
% Therefore, for each fixed power $t^N$, only finitely many terms in the series
% for $\Phi^{-1}$ contribute to its (i didnt understand this)).  The series is consequently a
% well-defined element of $\mathcal E$.  Multiplying formally gives
% \[
%   (\id_A+R)(\id_A-R+R^{\circ2}-R^{\circ3}+\cdots)=\id_A,
% \]
% and similarly on the other side.
% \end{proof}

% \begin{proposition}
% Suppose $\Phi_i\to\Phi$ in $\mathcal E$, and suppose each $\Phi_i$ has the
% form
% \[
%   \Phi_i=\id_A+O(t).
% \]
% Then
% \[
%   \Phi_i^{-1}\longrightarrow\Phi^{-1}.
% \]
% \end{proposition}

% \begin{proof}
% The limit has the form $\Phi=\id_A+O(t)$ and is therefore invertible by the
% preceding lemma.  We use the identity
% \[
%   \Phi_i^{-1}-\Phi^{-1}
%   =\Phi_i^{-1}\circ(\Phi-\Phi_i)\circ\Phi^{-1}.
% \](is this correct?)
% All factors other than $\Phi-\Phi_i$ have order $0$.  Hence
% \[
%   \ord(\Phi_i^{-1}-\Phi^{-1})
%   \geq\ord(\Phi_i-\Phi)\longrightarrow\infty.
% \]
% Thus $\Phi_i^{-1}\to\Phi^{-1}$.
% \end{proof}

\subsection{Continuity of the action on multiplication}

We now consider the action of formal endomorphisms on formal
multiplications. While $\mathcal E$ consists of formal endomorphisms of
$A$, a multiplication is a formal series of bilinear maps from
$A\otimes A$ to $A$. We therefore introduce
\[
\mathcal M
=
Hom_K(A\otimes A,A)[[t]].
\]

Thus,
\[
\mu\in\mathcal M,
\qquad
\phi,\psi,\psi'\in\mathcal E.
\]
The space $\mathcal M$ is equipped with the same $t$-adic metric and
order as above.

For
\[
\phi\in\mathcal E,
\qquad
\mu\in\mathcal M,
\qquad
\psi,\psi'\in\mathcal E,
\]
define
\[
T(\phi,\mu,\psi,\psi')
=
\phi\circ\mu\circ(\psi\otimes\psi').
\]
Equivalently, for $a,b\in A$,
\[
T(\phi,\mu,\psi,\psi')(a,b)
=
\phi\bigl(\mu(\psi(a),\psi'(b))\bigr).
\]

\begin{proposition}
Suppose
\[
\phi_i\to\phi,
\qquad
\mu_i\to\mu,
\qquad
\psi_i\to\psi,
\qquad
\psi_i'\to\psi'
\]
in $\mathcal E$, $\mathcal M$, $\mathcal E$, and $\mathcal E$,
respectively. Then
\[
T(\phi_i,\mu_i,\psi_i,\psi_i')
\longrightarrow
T(\phi,\mu,\psi,\psi').
\]
\end{proposition}

\begin{proof}
For convenience, write
\[
T_i
=
T(\phi_i,\mu_i,\psi_i,\psi_i')
\]
and
\[
T
=
T(\phi,\mu,\psi,\psi').
\]

We compare $T_i$ and $T$ by changing one factor at a time. We have the
decomposition
\[
T_i-T
=
T(\phi_i-\phi,\mu_i,\psi_i,\psi_i')+
T(\phi,\mu_i-\mu,\psi_i,\psi_i')+
T(\phi,\mu,\psi_i-\psi,\psi_i')+
T(\phi,\mu,\psi,\psi_i'-\psi').
\]

Indeed, expanding the right-hand side gives
\[(\phi_i-\phi)\mu_i(\psi_i\otimes\psi_i')
+\phi(\mu_i-\mu)(\psi_i\otimes\psi_i')
+\phi\mu\bigl((\psi_i-\psi)\otimes\psi_i'\bigr)
+\phi\mu\bigl(\psi\otimes(\psi_i'-\psi')\bigr),
\]
and the intermediate terms cancel, leaving
\[
\phi_i\mu_i(\psi_i\otimes\psi_i')
-
\phi\mu(\psi\otimes\psi').
\]

We now estimate the order of each term. By the order estimate for
composition and tensor products, we have
\[
\begin{aligned}
\ord\bigl(
T(\phi_i-\phi,\mu_i,\psi_i,\psi_i')
\bigr)
\geq{}&
\ord(\phi_i-\phi)
+\ord(\mu_i)
\\
&+
\ord(\psi_i)
+\ord(\psi_i').
\end{aligned}
\]

Similarly,
\[
\begin{aligned}
\ord\bigl(
T(\phi,\mu_i-\mu,\psi_i,\psi_i')
\bigr)
\geq{}&
\ord(\phi)
+\ord(\mu_i-\mu)
\\
&+
\ord(\psi_i)
+\ord(\psi_i'),
\end{aligned}
\]
while
\[
\begin{aligned}
\ord\bigl(
T(\phi,\mu,\psi_i-\psi,\psi_i')
\bigr)
\geq{}&
\ord(\phi)
+\ord(\mu)
\\
&+
\ord(\psi_i-\psi)
+\ord(\psi_i'),
\end{aligned}
\]
and
\[
\begin{aligned}
\ord\bigl(
T(\phi,\mu,\psi,\psi_i'-\psi')
\bigr)
\geq{}&
\ord(\phi)
+\ord(\mu)
\\
&+
\ord(\psi)
+\ord(\psi_i'-\psi').
\end{aligned}
\]

Since
\[
\phi_i\to\phi,\qquad
\mu_i\to\mu,\qquad
\psi_i\to\psi,\qquad
\psi_i'\to\psi',
\]
we have
\[
\ord(\phi_i-\phi)\longrightarrow+\infty,
\]
\[
\ord(\mu_i-\mu)\longrightarrow+\infty,
\]
\[
\ord(\psi_i-\psi)\longrightarrow+\infty,
\]
and
\[
\ord(\psi_i'-\psi')\longrightarrow+\infty.
\]

Moreover, all the remaining factors have nonnegative order, since they
are formal power series in nonnegative powers of $t$. Hence each of the
four terms in the decomposition has order tending to
$+\infty$.

Using
\[
\ord(f+g)
\geq
\min\{\ord(f),\ord(g)\},
\]
repeatedly, we obtain
\[
\ord(T_i-T)
\longrightarrow+\infty.
\]
Therefore,
\[
T(\phi_i,\mu_i,\psi_i,\psi_i')
\longrightarrow
T(\phi,\mu,\psi,\psi').
\]
\end{proof}
% \subsection{Continuity of the action on multiplication}

% For formal maps of the appropriate types define
% \[
%   T(\phi,\mu,\psi,\psi')
%   :=\phi\,\mu\,(\psi\otimes\psi').
% \]

% \begin{proposition}
% If
% \[
%   \phi_i\to\phi,\qquad
%   \mu_i\to\mu,\qquad
%   \psi_i\to\psi,\qquad
%   \psi_i'\to\psi',
% \]
% then
% \[
%   T(\phi_i,\mu_i,\psi_i,\psi_i')
%   \longrightarrow
%   T(\phi,\mu,\psi,\psi').
% \]
% \end{proposition}

% \begin{proof}
% Subtract the two expressions and insert intermediate terms.  The difference
% is a sum of terms in which exactly one of the four factors is replaced by the
% corresponding difference.  For example,
% \[
%   (\phi_i-\phi)\mu_i(\psi_i\otimes\psi_i')
% \]
% is of order tending to $+\infty$ because $\phi_i-\phi\to0$, while the other
% factors have nonnegative order.  The same argument applies to the remaining
% three terms.  Hence the total difference tends to zero.(This proof needs to be redone, more exposition)
% \end{proof}
\subsection{Infinite composition of elementary transformations}

The preceding results allow us to make sense of an infinite sequence of
elementary formal transformations. Suppose
\[
\phi_i={id}_A+q_i t^{m_i},
\qquad q_i\in {Hom}_K(A,A),
\]
where
\[
1\le m_0<m_1<m_2<\cdots.
\]
Define the partial compositions by
\[
\Phi_0={id}_A,\qquad
\Phi_i=\phi_{i-1}\circ\cdots\circ\phi_1\circ\phi_0
\quad(i\ge1).
\]
The strictly increasing positive exponents \(m_i\) ensure that each successive transformation affects only progressively higher powers of \(t\). The
following lemma makes this precise.

\begin{lemma}
For each $i\geq1$, there exist
\[
f_0,f_1,\ldots,f_{i-1}\in Hom_K(A,A)
\]
such that
\[
\Phi_i
=
\id_A
+
\sum_{r=0}^{i-1}f_r t^{m_r}
+
E_i,
\]
where
\[
\ord(E_i)>m_{i-1}.
\]
\end{lemma}

\begin{proof}
We proceed by induction on $i$.

For $i=1$, we have
\[
\Phi_1=\phi_0
=
\id_A+q_0t^{m_0}.
\]
Thus the required form holds with
\[
f_0=q_0
\]
and
\[
E_1=0.
\]

Now suppose that, for some $i\geq1$,
\[
\Phi_i
=
\id_A+\sum_{r=0}^{i-1}f_r t^{m_r}+E_i,
\]
where
\[
\ord(E_i)>m_{i-1}.
\]
Since
\[
\Phi_{i+1}=\phi_i\circ\Phi_i,
\]
we have
\[
\begin{aligned}
\Phi_{i+1}
&=
(\id_A+q_i t^{m_i})
\circ
\left(
\id_A+\sum_{r=0}^{i-1}f_r t^{m_r}+E_i
\right).
\end{aligned}
\]

Expanding the composition gives
\[
\Phi_{i+1}
=
\id_A
+
\sum_{r=0}^{i-1}f_r t^{m_r}
+
q_i t^{m_i}
+
E_{i+1},
\]
where $E_{i+1}$ contains all the remaining terms.

Indeed, every term arising from
\[
q_i t^{m_i}
\circ
f_r t^{m_r}
\]
has order
\[
m_i+m_r>m_i,
\]
and similarly
\[
\ord(q_i t^{m_i}\circ E_i)
>
m_i.
\]
Together with
\[
\ord(E_i)>m_{i-1}
\]
and
\[
m_i>m_{i-1},
\]
this gives
\[
\ord(E_{i+1})>m_i.
\]

Thus
\[
\Phi_{i+1}
=
\id_A
+
\sum_{r=0}^{i}f_r t^{m_r}
+
E_{i+1},
\]
with
\[
f_i=q_i
\]
and
\[
\ord(E_{i+1})>m_i.
\]
This completes the induction.
\end{proof}

\begin{lemma}
The sequence $(\Phi_i)$ is Cauchy in $\mathcal E$.
\end{lemma}

\begin{proof}
Let $N\geq0$ be given. Since
\[
m_0<m_1<m_2<\cdots
\]
is a strictly increasing sequence of nonnegative integers, we have
\[
m_i\longrightarrow+\infty.
\]
Choose $I$ sufficiently large so that
\[
m_{I-1}>N.
\]

For $i,j\geq I$, the preceding lemma gives
\[
\Phi_i
=
\id_A+\sum_{r=0}^{i-1}f_r t^{m_r}+E_i
\]
with
\[
\ord(E_i)>m_{i-1}\geq m_{I-1}>N,
\]
and similarly
\[
\Phi_j
=
\id_A+\sum_{r=0}^{j-1}f_r t^{m_r}+E_j
\]
with
\[
\ord(E_j)>m_{j-1}\geq m_{I-1}>N.
\]

Since $i,j\geq I$, the coefficients of $\Phi_i$ and $\Phi_j$ up to
degree $N$ are the same. Hence
\[
\Phi_i\equiv\Phi_j\pmod{t^{N+1}},
\]
and therefore
\[
\ord(\Phi_i-\Phi_j)>N.
\]

It follows that
\[
d(\Phi_i,\Phi_j)<2^{-N}.
\]
Since $N$ was arbitrary, $(\Phi_i)$ is a Cauchy sequence in
$\mathcal E$.
\end{proof}

\begin{proposition}
There exists a unique formal endomorphism
\[
\Phi\in\mathcal E
\]
such that
\[
\Phi_i\longrightarrow\Phi.
\]
Moreover,
\[
\Phi=\id_A+\sum_{r\geq0}f_r t^{m_r},
\]
and $\Phi$ is invertible with
\[
\Phi_i^{-1}\longrightarrow\Phi^{-1}.
\]
\end{proposition}

\begin{proof}
By the preceding proposition, $(\Phi_i)$ is a Cauchy sequence in
\[
\mathcal E=Hom_K(A,A)[[t]].
\]
The completeness result established earlier shows that, for each fixed
power $t^k$, the coefficient of $t^k$ in $\Phi_i$ is eventually constant.
Thus the stabilized coefficients define an element
\[
\Phi\in\mathcal E
\]
and
\[
\Phi_i\longrightarrow\Phi.
\]

More specifically, by the preceding lemma, for each $i$,
\[
\Phi_i
=
\id_A+\sum_{r=0}^{i-1}f_r t^{m_r}+E_i,
\qquad
\ord(E_i)>m_{i-1}.
\]
Hence, for each fixed $r$, the coefficient of $t^{m_r}$ is eventually
equal to $f_r$. Therefore the limit has the form
\[
\Phi
=
\id_A+\sum_{r\geq0}f_r t^{m_r}.
\]
The uniqueness of $\Phi$ follows from the fact that $d$ is a metric.

Since every $\Phi_i$ has the form
\[
\Phi_i=\id_A+O(t),
\]
the limit also has the form
\[
\Phi=\id_A+O(t).
\]
By the formal inverse lemma, $\Phi$ is therefore invertible.

Finally, by the continuity of inversion established in the preceding
subsection,
\[
\Phi_i^{-1}\longrightarrow\Phi^{-1}.
\]

\end{proof}

\section{Infinite elimination of coboundaries and Gerstenhabers Proposition 1}

We now return to the elimination process (\ref{elimination lemma}) for a formal deformation given in section \ref{formal deformation}. Let
\[
\mu^{(0)}
=
\mu+\sum_{k\geq1}F_k t^k
\]
be a one-parameter formal deformation of the multiplication $\mu$.

Suppose that at some stage of the elimination process we have
\[
\mu^{(r)}
=
\mu+t^{m_r}F_{m_r}+O(t^{m_r+1}),
\]
where $F_{m_r}$ is the first nonzero coefficient. Associativity implies that
\[
F_{m_r}\in Z^2(A,A).
\]

If
\[
[F_{m_r}]\neq0\in H^2(A,A),
\]
then the elimination process stops. The deformation has been transformed
into an equivalent deformation whose first nonvanishing coefficient
represents a nonzero cohomology class.

The remaining case is therefore when
\[
[F_{m_r}]=0.
\]
Then
\[
F_{m_r}=\delta q_r
\]
for some
\[
q_r\in Hom_K(A,A).
\]
Define the elementary formal transformation
\[
\phi_r=\id_A+q_r t^{m_r}.
\]
The elimination lemma gives
\[
\mu^{(r+1)}
=
\phi_r\mu^{(r)}
(\phi_r^{-1}\otimes\phi_r^{-1})
=
\mu+O(t^{m_r+1}).
\]
Thus the next nonzero term occurs at an order strictly larger than
$m_r$. Writing its order as $m_{r+1}$, we have
\[
m_{r+1}>m_r.
\]

Consequently, if the elimination process does not terminate, we obtain
\[
m_0<m_1<m_2<\cdots
\]
and a sequence of elementary transformations
\[
\phi_r=\id_A+q_r t^{m_r}.
\]

\subsection{The infinite case}

Assume that the elimination process continues indefinitely. Recall the
partial compositions
\[
\Phi_0=\id_A,
\qquad
\Phi_i
=
\phi_{i-1}\circ\cdots\circ\phi_1\circ\phi_0
\qquad (i\geq1).
\]

We have seen $(\Phi_i)$ is Cauchy and the completeness and coefficientwise
construction established there therefore give a formal endomorphism
\[
\Phi\in\mathcal E
\]
such that
\[
\Phi_i\longrightarrow\Phi.
\]
Moreover,
\[
\Phi_i^{-1}\longrightarrow\Phi^{-1}.
\]

We next relate the partial transformations to the successive
multiplications. By induction on $i$,
\[
\mu^{(i)}
=
\Phi_i\mu^{(0)}
(\Phi_i^{-1}\otimes\Phi_i^{-1}).
\]
Indeed, this is clear for $i=0$. If it holds for $i$, then
\[
\begin{aligned}
\mu^{(i+1)}
&=
\phi_i\mu^{(i)}
(\phi_i^{-1}\otimes\phi_i^{-1})\\
&=
\phi_i
\Phi_i\mu^{(0)}
(\Phi_i^{-1}\otimes\Phi_i^{-1})
(\phi_i^{-1}\otimes\phi_i^{-1})\\
&=
\Phi_{i+1}\mu^{(0)}
(\Phi_{i+1}^{-1}\otimes\Phi_{i+1}^{-1}),
\end{aligned}
\]
because
\[
\Phi_{i+1}=\phi_i\circ\Phi_i.
\]

On the other hand, the elimination construction gives
\[
\mu^{(i)}
=
\mu+O(t^{m_i}).
\]
Since
\[
m_i\longrightarrow+\infty,
\]
we have
\[
\ord(\mu^{(i)}-\mu)\longrightarrow+\infty.
\]
Hence
\[
\mu^{(i)}\longrightarrow\mu.
\]

Now apply the continuity of the action on multiplication established in
the preceding section. From
\[
\Phi_i\longrightarrow\Phi,
\qquad
\mu^{(0)}\longrightarrow\mu^{(0)},
\qquad
\Phi_i^{-1}\longrightarrow\Phi^{-1},
\]
we obtain
\[
\mu^{(i)}
=
\Phi_i\mu^{(0)}
(\Phi_i^{-1}\otimes\Phi_i^{-1})
\longrightarrow
\Phi\mu^{(0)}
(\Phi^{-1}\otimes\Phi^{-1}).
\]

But we have already shown that
\[
\mu^{(i)}\longrightarrow\mu.
\]
Since the metric space is Hausdorff, the limit is unique. Therefore
\[
\mu
=
\Phi\mu^{(0)}
(\Phi^{-1}\otimes\Phi^{-1}).
\]

Thus
\[
\mu^{(0)}
=
\Phi^{-1}\mu(\Phi\otimes\Phi),
\]
so the original deformation $\mu^{(0)}$ is equivalent to the original
multiplication $\mu$.

We have therefore proved that an infinite sequence of eliminations can
occur only when the deformation is equivalent to the undeformed
multiplication.

\subsection{Gerstenhaber's Proposition 1}

We can now state a formulation of Gerstenhaber's Proposition~1
adapted to the elimination procedure developed above.

\begin{proposition}[Gerstenhaber, Proposition 1]
Let
\[
\mu^{(0)}
=
\mu+\sum_{k\geq1}F_k t^k
\]
be a one-parameter formal deformation of the associative algebra
$(A,\mu)$. If $\mu^{(0)}$ is not equivalent to the original multiplication
$\mu$, then $\mu^{(0)}$ is equivalent to a deformation of the form
\[
\bar\mu
=
\mu+t^nF_n+t^{n+1}F_{n+1}+\cdots,
\]
where
\[
F_n\in Z^2(A,A)
\qquad\text{and}\qquad
[F_n]\neq0\in H^2(A,A).
\]
\end{proposition}

\begin{proof}
Let $m_0$ be the order of the first nonzero coefficient of
$\mu^{(0)}-\mu$. Associativity implies that
\[
F_{m_0}\in Z^2(A,A).
\]

If
\[
[F_{m_0}]\neq0,
\]
then the required form has already been obtained.

Suppose instead that
\[
[F_{m_0}]=0.
\]
Then
\[
F_{m_0}=\delta q_0
\]
for some $q_0\in Hom_K(A,A)$, and the elimination lemma gives an equivalent
deformation
\[
\mu^{(1)}
=
\mu+O(t^{m_0+1}).
\]
If the first nonzero coefficient of $\mu^{(1)}-\mu$ has a nonzero
cohomology class, the process stops and the required form is obtained.

Otherwise, the same argument can be repeated. Thus there are two
possibilities.

First, the process terminates after finitely many steps. At the final
stage, the first nonzero coefficient is a cocycle whose cohomology class
is nonzero, giving the required deformation.

Second, the process continues indefinitely. In that case the orders
satisfy
\[
m_0<m_1<m_2<\cdots,
\]
and the preceding infinite-elimination argument shows that the original
deformation $\mu^{(0)}$ is equivalent to the original multiplication
$\mu$.

Since, by hypothesis, $\mu^{(0)}$ is not equivalent to $\mu$, the second
case is impossible. Hence the elimination process must terminate after
finitely many steps, and the resulting equivalent deformation has the
form
\[
\bar\mu
=
\mu+t^nF_n+t^{n+1}F_{n+1}+\cdots,
\]
with
\[
F_n\in Z^2(A,A)
\qquad\text{and}\qquad
[F_n]\neq0\in H^2(A,A).
\]
\end{proof}

% \begin{remark}
% The infinite case is essential in this argument. The increasing sequence
% of elimination orders ensures that the corresponding elementary
% transformations converge to a single formal automorphism. At the same
% time, the transformed multiplications converge back to the original
% multiplication. The continuity of the action then shows that the
% infinite sequence of equivalences is itself represented by one formal
% equivalence.
% \end{remark}

\begin{corollary}
If
\[
H^2(A,A)=0,
\]
then every one-parameter formal deformation of $A$ is equivalent to the
original multiplication. Thus $A$ is rigid with respect to one-parameter
formal deformations.
\end{corollary}

\begin{proof}
Suppose that $\mu^{(0)}$ is a formal deformation which is not equivalent
to $\mu$. By Gerstenhaber's Proposition~1, it is equivalent to a
deformation whose first nonzero coefficient $F_n$ satisfies
\[
[F_n]\neq0\in H^2(A,A).
\]
This contradicts
\[
H^2(A,A)=0.
\]
Therefore every formal deformation is equivalent to the original
multiplication.
\end{proof}

\begin{remark}
The converse does not hold in general: nonvanishing of
$H^2(A,A)$ does not by itself imply the existence of a nontrivial formal
deformation. A nonzero cohomology class need not extend to a full
deformation.
\end{remark}

\end{document}